\documentclass[a4paper,11pt,reqno]{amsart}

\usepackage{amsmath}
\usepackage{amssymb}
\usepackage{xspace,epsfig}
\usepackage{mathptmx}
\usepackage{colortbl}

\usepackage{subfig}
\usepackage{color}

\usepackage{amsmath,amsfonts,enumerate,epsfig,amssymb,lscape}
\usepackage{amssymb}

\usepackage{amsmath,amsfonts,amssymb,mathrsfs,amsthm}

\numberwithin{equation}{section}
\usepackage{hyperref}

\newtheorem{Rem}{Remark}
\newtheorem{Def}{Definition}
\newtheorem{Lem}{Lemma}
\newtheorem{Theo}{Theorem}
\newtheorem{Prop}{Proposition}
\numberwithin{equation}{section}

\begin{document}

\title{Companion paper\\
 for\\
   Well posedness of general cross-diffusion systems}

\author[]{Catherine Choquet}
\address{$^a$ La Rochelle Univ., MIA, Avenue A. Einstein, F-17031, La Rochelle, France.
$^b$ CNRS EA 3165, France}
\email{cchoquet@univ-lr.fr}

%\author[C. Rosier]{Carole Rosier}
%\address{$^c$ Univ. Littoral C\^ote d'Opale, UR 2597, LMPA, Laboratoire de Math\'ematiques Pures et Appliqu\'ees Joseph Liouville, F-62100 Calais, France. 
%$^d$ CNRS FR 2956, France}
%\email{rosier@univ-littoral.fr}

%\author[L. Rosier]{Lionel Rosier}
%\address{$^c$ Univ. Littoral C\^ote d'Opale, UR 2597, LMPA, Laboratoire de Math\'ematiques Pures et Appliqu\'ees Joseph Liouville, F-62100 Calais, France.
%$^d$ CNRS FR 2956, France}
%\email{Lionel.Rosier@univ-littoral.fr}

\date{}

\begin{abstract}
The paper entitled ``Well posedness of general cross-diffusion systems'' \cite{CRR} is devoted to the mathematical analysis of the Cauchy problem for 
general cross-diffusion systems without any assumption about its entropic structure.
The absence of this type of hypothesis is strongly felt for two questions: the uniqueness of the solution, despite the  nonlinear coupling of the  highest derivatives terms, and the maximum principle.
The article  \cite{CRR} is therefore largely devoted to these two points.
The answers are  provided at the cost of certain assumptions or technicalities, mainly:
\begin{itemize}
\item  the {\it ratios} between the diffusion and cross-diffusion coefficients has to be drastically controlled for sufficiently enhancing the regularity of the solution, namely its gradient belongs to  the space $ L ^ 4 ((0,T)\times \Omega )$; the regularity is obtained by adapting the classical Meyer's to the nonlinear parabolic setting under consideration ;  
\item the  source terms have to  ensure the confinement of the solution.
\end{itemize}
The present ``companion'' paper aims at showing where more classical analysis tools fail to solve these questions and gives some additional clarifications.
\end{abstract}
%

%\vspace{0.3cm}
%\textbf{2010 Mathematics Subject Classification: 35K40, 93B05}  

\maketitle

\vspace{0.5cm}

%\begin{keyword}
%% keywords here, in the form: keyword \sep keyword
%Modeling of a seawater intrusion problem  \sep numerical simulations \sep finite element method \sep  a priori estimations.
%% MSC codes here, in the form: \MSC code \sep code
%% or \MSC[2008] code \sep code (2000 is the default)
\textbf{Keywords:}    cross-diffusion system; quasilinear parabolic equations; uniqueness in the small; boundedness. 
%\end{keyword}

%%%%%%%%%%%%%%%%%%%%%%%%%%%%%%%%%%%%%%%%%%%%%%%%%%%%%%%%%%%%%%%%%%%%%%%%%%%%%%%%%%%%%%%%%%%%%%%%%%%%%%%%%%%%%%%%%%%%%%%%%%%%%%%

%% The amsthm package provides extended theorem environments
%% \usepackage{amsthm}

\section{Introduction}
\label{s1}

In what follows, excerpts from the article will be written in blue.
Some notations are recalled in the present section.
The second section is devoted to the uniqueness result and some points related to the maximum principle are presented in the third section.

\bigskip

\textcolor{blue}{
We consider an open bounded domain $\Omega$ of $\mathbb{R}^N$, $N \in \mathbb{N}^*$, $N \le 3$.
The boundary of $\Omega$, assumed to be of class  ${\mathcal C}^1$, is denoted by  $\Gamma$. 
The time interval of interest is $(0,T)$, $T$ being any positive real number. 
Set $\Omega _T := (0,T) \times \Omega$.
All the results and all the computations of \cite{CRR} are done for a particular class of cross-diffusion systems, the one {classically modeling the dispersal of two interacting biological species.}
Indeed, it is one of the less cumbersome systems containing all the difficulties inherent to the analysis of a strongly coupled cross-diffusion:
\begin{equation}   
\label{model}
%\displaystyle 
\partial_t u_i -\nabla \cdot \big(\delta_i  \nabla u_i + u_i \sum_{j=1}^m K_{i,j}  \nabla u_j \big) 
=  Q_i(u) \  \textrm{ in }\Omega_T, \textrm{ for }    i= 1,...,m.
\end{equation}
It is completed by the following boundary and initial conditions, for $i=1,..., m$:
$$
 u_i= u_{i,D} \ \mathrm{in} \ (0,T)\times \Gamma , 
\quad
 u_i(0,x)=u_{i} ^0 (x) \ \mathrm{in} \ \Omega .
$$  
For any $1\le i,j\le m$, the  tensor $K_{i,j}$ is assumed to be bounded and uniformly elliptic. More precisely, there exist two positive real numbers, $0< K_{i,j}^- \le K_{i,j}^+$, such that
\begin{equation}
       0<K_{i,j}^-|\xi|^2\leq K_{i,j} \xi \cdot \xi =\sum_{k,l=1}^N (K_{i,j})_{kl} \xi_k\xi_l\leq K_{i,j}^+|\xi|^2 ,\ \  \forall \xi\in\mathbb{R}^N\setminus \{ 0\} .
\label{AssK}       
\end{equation}
We consider the fully non-degenerate setting
\begin{equation}
\delta_i >0 \quad 1\le i \le m,
\label{delta}
\end{equation}
thus prohibiting the full exploitation of entropy methods. 
}

The previous sentence deserves some attention. 
First of all, it is interesting from a pedagogical point of view: thinking that a parabolic problem is easier to analyze than a degenerate parabolic problem is indeed sometimes misleading.
Such an {\it a priori} is tempting because of the regularity result, in $L^2(0,T;H^1(\Omega))$, which is usually induced by the parabolic structure. 
But it is still necessary to be able to demonstrate that a solution exists for it to inherit this regularity!
Moreover, as explained in the Introduction of \cite{CRR}, losing the entropic structure also makes us lose one of the usual methods to prove a maximum principle for the solution.
Let's add to the confusion.
As already mentioned, the system considered here may be viewed as a model for the dispersal of two interacting biological species. 
Its degenerate setting is partly considered  by Carrillo et al. in \cite{Car}. 
They write ``The main mathematical difficulty here arises from the cross-diffusion term allowing for segregation fronts\footnote{Notice that such segregation fronts do not make sense in many physical situations, thus the importance of considering the non-degenerate setting. }  to form in the solutions.[...] These remarkable results have severe consequences, initially smooth solutions lose their regularity when both densities meet each other. In fact, they become discontinuous at the contact interface immediately.'' 
Why do not they share our analysis of the difficulty? 
Because they have to face a kind of `ultimate maximum principle', namely a segregative result: if one of the unknowns reaches a given maximal value, the other one vanishes. 
Here, on the contrary, a simple boundedness result requires far from obvious considerations.

\medskip
\textcolor{blue}{
Let us now introduce some elements for the functional setting used in the present paper. 
For the sake of brevity we shall write $H^1(\Omega )=W^{1,2}(\Omega )$ and 
\begin{eqnarray*}
V=H_0^1(\Omega ),\ V'=H^{-1}(\Omega ),\ H=L^2(\Omega ).
\end{eqnarray*}
The embeddings  $V\subset H=H'\subset V'$ are dense and compact. 
For any $T>0$, let $W(0,T)$ denote the space 
$$
W(0,T):=
\bigl\{ \omega \in L^2(0,T;V ),\ \partial_t \omega  \in L^2(0,T;V ')\bigr\}
$$
endowed with the Hilbertian norm 
$ \Vert \omega\Vert _{W(0,T)}^2 =
 \Vert \omega  \Vert ^2_{L^2(0,T;V )}+
\Vert \partial_t \omega \Vert ^2_{ L^2(0,T;V ')} $.
We assume that there exists a lifting of each boundary function $u_{i,D}$, still denoted the same for convenience, belonging to the space $L^2(0,T;H^1(\Omega)) \cap H^1(0,T;(H^1(\Omega))')$. Due to the smoothness of $\Gamma$, such a result is ensured if $u_{i,D} \in L^2(0,T;H^{1/2}(\Gamma)) \cap H^1(0,T;H^{-1/2}(\Gamma)))$ (see \cite{LiMa}).
The initial data $u_{i}^0$  are assumed to be in  $H$, the source terms $Q_i(v)$  to be in $L^2(\Omega_T)$ for any $v \in (W(0,T))^m$, $ 1\le i \le m$.  
}

%%%%%%%%%%%%%%%%%%%%%%%%%%%%%%%%%%%%%%%%%%%%%%%%%%

\medskip

%%%%%%%%%%%%%%%%%%%%%
\section{Uniqueness}
\label{s3}
%%%%%%%%%%%%%%%%%%%%%

Here, as in the article, we postpone the difficulty related to the establishment of a  maximum principle to later and we  begin by considering the problem with some bounded nonlinearities.
To this aim, we introduce , for $\ell>0$, the truncating function $T_\ell$ defined by
$$
T_\ell(u)=  \max\bigl\{ 0, \min \{ u,\ell \} \bigr\} . $$
%and we extend the function $T_\ell$ constantly in to $\mathbb{R}$.
We then consider the following problem: for $ i=1, ... ,m$, 
\begin{eqnarray}
&& \partial_t u_i -\nabla \cdot \big(\delta_i \, \nabla u_i\,+ \, T_\ell(u_i ) \sum_{j=1}^m K_{i,j}  \nabla u_j \big) =  Q_i(u)
\quad  \textrm{in }\Omega_T,
\label{problem} \\ 
&& u_i= u_{i,D} \ \mathrm{in} \ (0,T)\times \Gamma , 
%\label{boundary}
\quad
 u_i(0,x)=u_{i} ^0 (x)\  \mathrm{in} \ \Omega . 
\label{initial}
\end{eqnarray}  
The initial and boundary conditions are supposed to satisfy the compatibility conditions
  \begin{equation}
  \label{FF1}
  u_{i} ^0(x)=u_{i,D}(0,x), \quad x \in \Gamma , \quad 1\le i \le m .
  \end{equation}
For the sake of simplicity, we set $m=2$.  
The following existence result is proved in \cite{CRR}.

\vspace{0.05cm}  
\textcolor{blue}{
\begin{Theo} \label{Theo1}
     Assume  that the tensor $K$ satisfies:
     \begin{equation}
      \dfrac{(K_{1,2}^+)^2}{K_{1,1}^-} < \dfrac{4\delta _2}{\ell},\, \, 
      \dfrac{(K_{2,1}^+)^2}{K_{2,2}^-} < \dfrac{4\delta _1}{\ell} .
        \label{AA1bis}
        \end{equation}
 Then for any $ T > 0 $,  the problem \eqref{problem}--\eqref{initial} admits a weak 
solution  $(u_i)_{i=1,2} \in (W(0,T))^2$.
Furthermore,  if almost everywhere  in $\Omega_T$, $0 \le u^0_i $, $0 \le u_{i,D}$ and $Q_i(v) \ge 0$ if $v_i \le0$, the following relation holds true
$$ 0\leq u_i(t,x) \quad \textrm{for a.e. } x \in \Omega , \textrm{ for all } t \in (0,T), \  i=1,2 .$$
  \end{Theo}
}
 
   \smallskip

Proving a uniqueness result for a cross-diffusive system is always a tricky problem. 
In \cite{CRR}, the results are founded on an additional regularity result, namely  a Meyer's type property allowing to upgrade the regularity of any solution of the cross-diffusive problem from $L^2(H^1)$ to $L^4(W^{1,4})$.
Forcing the regularity of the solution in this way could seem unnatural since the typical Meyer's  result is an upgrading from $L^2(;H^1)$ to $L^s(W^{1,s})$, for some $s=2+\epsilon$, where $\epsilon >0$ could be {\it a priori} very small. 
We thus discuss  this result in the following subsection. 
The second subsection presents the kind of uniqueness result we can prove without forcing the regularity.

\smallskip

%%%%%%%%%%%%%%%%%%%%%
\subsection{Enhanced regularity result}
%%%%%%%%%%%%%%%%%%%%%

We begin by a parabolic extension of the Meyers regularity theorem \cite{Meyers}. 
Once again, we introduce some notations.
\textcolor{blue}{
Let $X_p=L^p(0,T;W_0^{1,p}(\Omega))$, $p \ge 2$, endowed with the norm 
\[
\Bigl(\int_0^T||v(t)||_{W_0^{1,p}(\Omega)}^p dt\Bigr)^{1/p} :=  ||\nabla v||_{L^p(\Omega_T)^N} .
\]
 The space $Y_p=L^p(0,T;W^{-1,p}(\Omega))$ is endowed with the norm $ ||f||_{Y_p}=\inf_{\textrm{div} _x g=f}||g||_{(L^p(\Omega_T))^N}$.
Given $F\in Y_p$,  there is a unique solution $u\in X_p$ of the following initial boundary value problem 
\begin{eqnarray*}
&& \partial_t u - \Delta u  =F\mbox{ in } \Omega_T,\quad
 u=0 \textrm{ on } (0,T)\times \Gamma, \quad    
u(0,x)=0 \textrm{ in } \Omega .
\end{eqnarray*}
We set $\Lambda^{-1} = \partial_t - \Delta$, so that  $u=\Lambda (F)$. 
Let $g$ be defined by
$${g}(p): =||\Lambda ||_{\mathcal{L}(Y_p;X_p)}.$$
 It is well-known  that  ${g}(2)=1$. 
Now, let $A\in (L^\infty(\Omega ))^{N\times N}$ be such that there exists $\alpha>0$ satisfying
\begin{eqnarray*}
\sum_{i,j=1}^N A_{i,j}(x)\xi_i\xi_j \geq \alpha |\xi|^2 \ \textrm{ for  a.e. } x\in\Omega \textrm{ and for all } \xi\in\mathbb{R}^N.
\end{eqnarray*}
We set $\beta :=\max_{1\leq i,j\leq n} ||A_{i,j}||_{L^\infty(\Omega)}$ and $\mathcal{A}u = -\sum_{i,j=1}^N \partial_{x_i} \bigl( A_{i,j} \partial_{x_j} u \bigr)$.
We  state the following Lemma ({\it cf} \cite{Bensoussan} and Appendix in \cite{CRR}).
%%%%%%%%%%%%%%%%%%%%%%%%%%%%%%%%%%%%%%%%%%%%%%%%%%
\begin{Lem}\label{Lem2}
%Assume that $\Gamma$ is smooth. 
Let $f\in L^2(0,T;V')$, $u^0\in H$ and $u \in L^2(0,T; V)$ be the solution of 
\begin{equation}
\partial_t u +\mathcal{A}u=f \mbox{ in } \Omega_T, \quad
u(0)=u^0 \mbox{ in } \Omega .
\label{PbLions}
\end{equation}
There exists $r>2$, depending on $\alpha,\beta$ and $\Omega$, such that if 
$u^0\in W_0^{1,r}(\Omega)$ and 
$ f\in Y_r$, 
then $u\in X_r$. Furthermore, the following estimate holds true
\begin{eqnarray}
||u||_{X_r}\leq {C}(\alpha,\beta,r)(||f||_{{Y_r}}+{\beta T^{1/r}}||u^0||_{W_0^{1,r}(\Omega)}),
\label{5.17}
\end{eqnarray}
where the  constant ${C}(\alpha,\beta,r)>0$ depends on $\Omega$, $\alpha$, $\beta$ and $r$ (but not on $T$) as follows: 
\begin{equation}
{C}({\alpha},{\beta},r)\leq \frac{ {g}(r)}{(1-{k}(r)) \,({\beta}+{c})} , \quad k(r)={g}(r)(1-{\mu}+{\nu})
\label{5.18} 
\end{equation}
where $\mu=(\alpha+c)/(\beta +c)$, $\nu=(\beta^2+c^2)^{1/2}/(\beta+c)$ and $c$ is any real number such that $c>(\beta^2-\alpha^2)/2\alpha$. 
If, moreover, $A$ is symmetric, the estimate \eqref{5.18} holds true with $\mu=\alpha/\beta$ and $\nu=c=0$.
\end{Lem}
%%%%%%%%%%%%%%%%%%%%%%%%%%%%%%%%%%%%%%%%%%%%%%%%%%
}

The latter lemma is actually a Meyers type result. 
Indeed, we have $1-\mu+\nu <1$. 
According to the Riesz-Thorin's theorem, the function $g$ is bounded by a continuous function $\rho$ such that $\rho(2)=g(2)=1$. 
It ensures that, if $s$ is close enough to 2, $k(s)<1$ thus the invertibility of the operator $\partial_t + \mathcal{A}$ from $X_s$ to $Y_s$.
The additional information here is a criterion, expressed with regard to the norm of the inverse of the Heat operator $\Lambda$, which basically states how close to $\Lambda$ the operator $\partial_t + \mathcal{A}$ has to be for ensuring its invertibility.

In the spirit of Meyers' regularity result, the existence result of a solution of  \eqref{problem}--\eqref{initial},  Proposition 1 in \cite{CRR}, could thus be rewritten as follows.
\begin{Prop}
\label{prop1}
Set ${\alpha_i}=\delta_i$, ${\beta_i}=\delta_i+{ \ell K_{i,i }^+}$,  ${\mu_i}=(\alpha_i+c_i)/(\beta_i+c_i)$ and ${\nu_i}^2=(\beta_i^2+c_i^2)/(\beta_i+c_i)^2$  for $i=1,2$.
Let $c_i=0$  if $K_{i,i}$ is  symmetric and  $c_i > (\beta_i^2-\alpha_i^2)/2\alpha_i$ if not. 
Let $(u_1,u_2)$ be a solution of Problem \eqref{problem}--\eqref{initial}. 
Assume that $(\ell , \delta_1, \delta_2)$ and the tensor $K$ satisfy 
\begin{eqnarray}
K_{i,j}^+ <  \frac{(\beta_i^*+c_i^*)  (  \mu_i^* - \nu_i^*) }{2\ell},   \quad i=1,2,\, i\neq j .
\label{CondReg}
\end{eqnarray}
Then, there exists some $s>2$ such that, if   $(u_{1}^0,u_{2}^0) \in (W^{1,s}(\Omega))^2$, then  $\nabla u_1$ and $\nabla u_2$ belong to  $(L^s(\Omega_T))^N$.
\end{Prop}

The reader may wonder what kind of uniqueness result may be obtained from the latter natural enhancement\footnote{that is without forcing the regularity to reach $s=4$}.
An answer is given in the following subsection.

\smallskip

%%%%%%%%%%%%%%%%%%%%%
\subsection{Uniqueness in the small result}
%%%%%%%%%%%%%%%%%%%%%

``Uniqueness in the small'' entitles Section 4.2 of the monograph \cite{Lady} co-authored by Olga A. Ladyzhenskaya and Nina U. Uralceva.
This work is especially renowned for providing a complete uniqueness analysis of quasilinear pde's in the scalar case. 
In the present setting, the following approach does not seem too presumptuous: 
we could transfer the uniqueness in the small result from p.257 of \cite{Lady} to the case of our system just following the same ideas and adding more and more restrictions when it becomes necessary.
Bear in mind that a (nonlocal) uniqueness result is provided in \cite{CRR} in the case $\Omega \subset \mathbb{R}^2$ as soon as a first restriction on the {\it ratios} between cross-diffusive and diffusive parameters ensures the additional regularity in $L^4(0,T;W^{1,4}(\Omega))$ of the solution and as another technical restriction is assumed. 
We now aim at checking if a local uniqueness result could be reached with weaker assumptions.
The computations are detailed in the following lines.
Notice that they also shed light on the assumptions that could lead to a result of overall uniqueness when $N>2$.

%\medskip

\subsubsection{Preliminary computations}\label{prel}

Assume that   $(u_1,u_2)$ and $(\bar u_1,\bar u_2)$  are two weak solutions of \eqref{problem}.
Then the functions 
$v_i:=u_i-\bar u_i\in W(0, T)$, $i=1,2$,  weakly solve the following system in $\Omega_T$: 
\begin{eqnarray*}
&&\partial_{t} v_1-\nabla\cdot\big( (\delta _1 +K_{1,1} T_\ell(u_1))\nabla v_1 \big)-\nabla\cdot  \big(K_{1,1} (T_\ell(u_1)-T_\ell(\bar u_1)\nabla \bar u_1\big) \\
&& \qquad \quad
-\nabla\cdot( K_{1,2} T_\ell(u_1) \nabla v_2 )-\nabla\cdot \big(K_{1,2} (T_\ell(u_1)-T_\ell(\bar u_1))\nabla \bar u_2\big)=0 ,
\\
&&\partial_{t} v_2 -\nabla\cdot\big( (\delta _2 +K_{2,2} T_\ell(u_2))\nabla v_2 \big)-\nabla\cdot \big(K_{2,2} (T_\ell(u_2)-T_\ell(\bar u_2))\nabla \bar u_2\big) \\
&& \qquad \quad
-\nabla\cdot( K_{2,1} T_\ell(u_2) \nabla v_1 )-\nabla\cdot \big(K_{2,1} (T_\ell(u_2)-T_\ell(\bar u_2))\nabla \bar u_1\big)=0 .
\end{eqnarray*}

Assume that $(u_1,u_2)$ and $(\bar u_1,\bar u_2)$ coincide a.e. $t \in (0,T)$ on the boundary $\partial K_\rho$ of a given open sphere $K_\rho \subset \Omega$ of radius $\rho$.
Then $v_1$ and $v_2$ satisfy homogeneous Dirichlet boundary conditions on $\partial K_\rho$.
We multiply the equations by, respectively, $v_1$ and $v_2$ and we integrate over $(0,t) \times K_\rho$ with $0<t\le T$. 
Using the fact that $v_1(0,.)=v_2(0,.)=0$ a.e. in $\Omega$ and the coercivity property of $K_{i,i}$, we get after summing up  the two equations:
\begin{eqnarray}
&&\frac{1}{2}\int_{K_\rho} \bigl( |v_1|^2(t,x) +|v_2|^2(t,x)\bigr)%\, dx 
%\\
%&&\quad 
+\int_0^t\int_{K_\rho} \bigl( (\delta_1+K_{1,1}^- T_\ell(u_1))  |\nabla v_1|^2
+(\delta _2 +K_{2,2}^- T_\ell(u_2))   |\nabla v_2|^2 \bigr) % \, dxds 
\nonumber \\ 
&&\quad 
+\int_0^t\int_{K_\rho} \bigl( T_\ell(u_1)-T_\ell(\bar u_1)\bigr)\big(K_{1,1}\nabla \bar u_1+K_{1,2} \nabla \bar u_2\big) \cdot\nabla v_1% \, dxds 
%\\
%&&\quad 
+\int_0^t \int_{K_\rho}   \big(K_{1,2} T_\ell(u_1)+K_{2,1} T_\ell(u_2)\big) \nabla v_1  \cdot\nabla v_2% \,  dxds 
\nonumber \\
&&\quad  
+ \int_0^t\int_{K_\rho}   \bigl(T_\ell(u_2)-T_\ell(\bar u_2)\bigr)\big(K_{2,1}\nabla \bar u_1+K_{2,2} \nabla \bar u_2\big) \cdot\nabla v_2% \, dxds 
\leq 0 .
\label{0}
\end{eqnarray}
Using the Cauchy-Schwarz and Young inequalities, we get for any arbitrary $\epsilon_{i+2} >0$, $i=1,2$: 
\begin{eqnarray*}
&& \Bigl|\int_0^t \int_{K_\rho} K_{i,-i} T_\ell(u_i) \nabla v_i  \cdot\nabla v_{-i} \Bigr|  
\leq \ell^{1/2}K_{i,-i}^+ \Bigl(\int_0^t \int_{K_\rho}|\nabla v_{-i} |^2\Bigr)^{1/2}\Bigl(\int_0^t\int_{K_\rho}T_\ell(u_i) |\nabla v_i |^2\Bigr)^{1/2} 
\\
&& \qquad \qquad \quad
\leq \frac{\ell (K_{i,-i}^+)^2 }{4\epsilon_{i+2}}(\int_0^t\int_{K_\rho}|\nabla  v_{-i}|^2)
+ \epsilon_{i+2}(\int_0^t\int_{K_\rho}T_\ell(u_i^*) |\nabla v_i|^2) .
\end{eqnarray*}
These terms may be treaten as in \cite{CRR} provided that $\ell (K_{i,-i}^+)^2 / K_{i,i}^-$ is sufficiently small with regard to $\delta_{-i}$.
We will therefore no longer pay attention to these terms.

By the definition of $T_\ell$ and since $u_i, \bar u_i \ge 0$,  we have that $T_\ell(u_i) \geq 0$ and   
 $ \vert T_{\ell}(u_i) -T_{\ell}(\bar u_i)\vert \leq \vert u_i - \bar u_i \vert = \vert v_i \vert$. 
 For notational convenience, let  $K_{i,+}=\max _{j=1,2}\vert K_{i,j}^+\vert$,  $i=1,2$.
We have 
\begin{eqnarray*}
&&\left\vert \int_0^t\int_{K_\rho}   (T_\ell(u_i)-T_\ell(\bar u_i))\big(K_{i,i}\nabla \bar u_i+K_{i,-i} \nabla \bar u_{-i}\big) \cdot\nabla v_i \, dxds \right \vert 
%\\
%&&\qquad 
\leq \int_0^t \int_{K_\rho}  K_{i,+} |v_i| \big(  | \nabla \bar u_i| + |\nabla \bar u_{-i}| \big) | \nabla v_i| \, dxds.
\end{eqnarray*}
All the difficulty induced by the cross-diffusive structure lies in the estimate of the latter integral, in the form
\begin{equation}
 I_{ij} = \int_0^t \int_{K_\rho}  K_{i,+} |v_i| \, | \nabla  u_j| \, | \nabla v_i| \, dxds, \  i,j=1,2.
 \label{defI}
 \end{equation}
According to  the Cauchy-Schwarz and Young inequalities, we have
$$ I_{ij} \le \epsilon_i \int_0^t\int_{K_\rho} \vert \nabla v_i\vert^2 + \frac{K_{i,+}}{4\epsilon_i}\int_0^t\int_{K_\rho} \vert \nabla u_j \vert^2 \vert v_i \vert^2$$
for any $\epsilon_i >0$ and {\bf we now focus on }
$$ J = \int_0^t\int_{K_\rho} \vert \nabla u_j \vert^2 \vert v_i \vert^2, \  i,j=1,2.$$

\medskip

\noindent{\bf Step 1: Estimate of  $\int_0^t\int_{K_\rho} \vert \nabla u_i \vert^2 $.}\\
Let $\varphi \in \mathcal{C}^\infty(\bar \Omega)$ such that $\varphi_{\vert K_\rho}=1$ and  $\varphi_{\vert ^c\!K_{\rho_1}}=0$ with $K_\rho \subset K_{\rho_1} \subset \Omega$, $\rho_1>\rho$.
We use the test function $u_i\varphi^2$, $i=1,2$, in the variational formulation of  \eqref{problem}.
We get for $i=1,2$:
\begin{eqnarray}
&& \int_0^t\int_\Omega \partial_t u_i u_i \varphi^2 
+\int_0^t \int_\Omega (\delta_i + K_{i,i}^- T_\ell(u_i)) \varphi^2 \vert \nabla u_i \vert^2
\nonumber \\
&& \qquad
+ \int_0^t \int_\Omega K_{i,j} T_\ell(u_i) \varphi^2 \nabla u_j \cdot \nabla u_i 
+ 2 \int_0^t \int_\Omega  (\delta_i + K_{i,i} T_\ell(u_i))u_i \varphi \nabla u_i \cdot \nabla \varphi
\nonumber \\
&& \qquad
+ 2 \int_0^t \int_\Omega K_{i,j} T_\ell(u_i) u_i \varphi \nabla u_j \cdot \nabla \varphi 
\le \int_0^t\int_\Omega Q u_i \varphi^2.
\label{1}
\end{eqnarray}
The term $ \int_0^t \int_\Omega K_{i,j} T_\ell(u_i) \varphi^2 \nabla u_j \cdot \nabla u_i $ may be controlled as in the existence proof of \cite{CRR}:
\begin{equation}
 \Bigl\vert  \int_0^t \int_\Omega K_{i,j} T_\ell(u_i) \varphi^2 \nabla u_j \cdot \nabla u_i \Bigr\vert
\le  \int_0^t \int_\Omega K_{i,i}^- T_\ell(u_i)\varphi^2 \vert \nabla u_i \vert^2 + \frac{\ell {K_{i,-i}^+}^2}{4K_{i,i}^-}  \int_0^t \int_\Omega  \varphi^2\vert \nabla u_{-i}\vert^2.
\label{3bis}
\end{equation}
We thus have to estimate the following quantities:
\begin{eqnarray*}
&& I_1 =  \int_0^t \int_\Omega \partial_t u_i u_i \varphi^2 = \frac{1}{2}  \int_0^t \frac{d}{dt} \int_\Omega u_i^2 \varphi^2,
\\
&& I_2 = 2 \int_0^t \int_\Omega(\delta_i + K_{i,i} T_\ell(u_i))u_i \varphi \nabla u_i \cdot \nabla \varphi,
\\
&& I_3 =  2 \int_0^t \int_\Omega K_{i,j} T_\ell(u_i) u_i \varphi \nabla u_j \cdot \nabla \varphi ,
\\
&& I_4 = \int_0^t\int_\Omega Q u_i \varphi^2 .
\end{eqnarray*}

 For $I_2$ and $I_3$, we have to estimate
 $$\int_0^t \int_\Omega u_i \vert \nabla u_j\vert \vert \nabla \varphi\vert = \int_0^t \int_{K_{\rho_1}} u_i \vert \nabla u_j\vert \vert \nabla \varphi\vert , \  j=i,-i.$$
 Set $\rho_1 = 2\rho$.
 We can choose $\varphi$ such that $\max_\Omega \vert \nabla \varphi\vert \le C/\rho$.
 Then
 \begin{eqnarray}
 &&  \int_0^t \int_{K_{2\rho}} u_i \vert \nabla u_j\vert \vert \nabla \varphi\vert 
 \le \frac{C}{\rho} \int_0^t \int_{K_{2\rho}} u_i \vert \nabla u_j\vert 
 \le \frac{C}{\rho} 
 \Bigl(  \int_0^t \int_{K_{2\rho}}  \vert \nabla u_j\vert^s \Bigr)^{1/s} 
 \Bigl( \int_0^t \int_{K_{2\rho}} u_i^{s/(s-1)} \Bigr)^{(s-1)/s}
 \nonumber \\
 && \qquad
 \le
 \frac{C}{\rho} C_s \Bigl( \int_0^t \Vert u_i \Vert_{L^\infty(\Omega)}^{s/(s-1)} \bigl( \int_{K_{2\rho}} dx\bigr) dt \Bigr)^{(s-1)/s}
 \le \frac{C}{\rho} C_s^2 (2\rho)^{N(s-1)/s} t^{(s-2)/s}
 \nonumber \\
 && \qquad
 = C C_s^2 t^{(s-2)/s} \rho^{((N-1)s-N)/s}
 \label{3}
 \end{eqnarray}
if we assume:
\begin{enumerate}
\item[(i)] a restriction on the {\it ratios} between cross-diffusive and diffusive parameters ensures the additional regularity in $L^s(0,T;W^{1,s}(\Omega))$, $s>2$,  of the solution and
$$ \Vert u_i \Vert_{L^s(0,T;W^{1,s}(\Omega))} \le C_s;$$
\item[(ii)] either the real number $s$ is large enough for ensuring the Sobolev injection $W^{1,s}(\Omega) \subset L^\infty(\Omega)$, that is 
$$ s>N.$$
\end{enumerate}
Notice that the estimate \eqref{3} may be replaced by the following
\begin{equation}
 \int_0^t \int_{K_{2\rho}} u_i \vert \nabla u_j\vert \vert \nabla \varphi\vert \le C C_s t^{(s-1)/s} \rho^{((N-1)s-N)/s}
 \label{3born}
\end{equation}
if we assume:
\begin{enumerate}
\item[(iii)] we deal with bounded solutions $u_i$.
\end{enumerate}

Assume for $I_1$ that the initial conditions in \eqref{problem} are such that
\begin{equation}
 \int_{K_{2\rho}} \vert u_{i}^0 \vert^2 \le C \rho^{((N-1)s-N)/s}.
 \label{2}
 \end{equation}
 
 Finally, the quantity $I_4$ may be controlled by $I_1$ thanks to the Gronwall lemma.
 Nevertheless, as we aim also deal with the time independent case, we provide another estimate. 
 We write for instance
$$ \vert I_4 \vert =\Bigl\vert  \int_0^t\int_\Omega Q u_i \varphi^2 \Bigr\vert  
 \le \Vert u_i  \Vert_{L^s(\Omega_t)} \Vert Q  \Vert_{L^{s/(s-1)}(K_{2\rho} \times (0,t))} 
 \le C_s \Vert Q  \Vert_{L^{s/(s-1)}(K_{2\rho} \times (0,t))} $$
 and we assume that $Q$ satisfies
 \begin{equation}
 \Vert Q  \Vert_{L^{s/(s-1)}(K_{2\rho} \times (0,t))}  \le C \rho^{((N-1)s-N)/s} .
 \label{2b}
 \end{equation}

We infer from \eqref{3bis}-\eqref{2b} in the sum of \eqref{1}$_i$, $i=1,2$, that
$$ \int_{\Omega_t } \varphi^2 \vert \nabla u_i \vert^2 \le C \bigl(\delta_i- \frac{\ell {K_{i,-i}^+}^2}{4K_{i,i}^-} \bigr)^{-1}C\bigl(C_s^2 t^{(s-2)/s} ,\sum_{j=1}^2\Vert u_j^0\Vert_{L^2}^2, \Vert Q  \Vert_{L^{s/(s-1)}(K_{2\rho} \times (0,t))}\bigr)\rho^{((N-1)s-N)/s}$$
and thus, in view of the definition of $\varphi$:
\begin{eqnarray}
 \int_0^t \int_{K_\rho} \vert \nabla u_i \vert^2 
&\le& C \bigl(\delta_i- \frac{\ell {K_{i,-i}^+}^2}{4K_{i,i}^-} \bigr)^{-1}C\bigl(C_s^2 t^{(s-2)/s} ,\sum_{j=1}^2\Vert u_j^0\Vert_{L^2(K_{2\rho})}^2, \Vert Q  \Vert_{L^{s/(s-1)}(K_{2\rho} \times (0,t))}\bigr) 
\nonumber \\
&&
\times \rho^{((N-1)s-N)/s}.
\label{4}
\end{eqnarray}
Notice that $\delta_i- \ell {K_{i,-i}^+}^2/(4K_{i,i}^-)>0$ in view of the assumption made for ensuring the existence of a solution for \eqref{problem} (see Theorem \ref{Theo1}).

\medskip

\noindent{\bf Step 2: Auxiliary results for turning back to $J$.}\\
We first mention Lemma 4.3 page 59 in \cite{Lady} (and its corolary) : if $m\ge 0$, if $\alpha>0$, if $\int_{K_\rho} \vert v \vert \le C \rho^{m+\alpha}$ then $\int_{K_\rho} \vert x-y \vert^{-m-\alpha/2} \vert v(x)\vert \, dx \le C_1(\alpha,m,C,\mbox{diam}(\Omega))\rho^{\alpha/2}$ for any $y \in K_\rho$. 
We have denoted by $\mbox{diam}(\Omega)$ the diameter of $\Omega$.
Set $v=\int_0^t \vert \nabla u_i \vert^2$. 
Set $m=N-2$ and $\alpha=(s-N)/s$. 
Assume $s>N$ so that $\alpha >0$.
We have $m+\alpha=((N-1)s-N)/s$.
We thus infer from \eqref{4} and Fubini's theorem that
\begin{eqnarray}
\int_0^t \int_{K_\rho} \vert x-y \vert^{-N+2-\alpha/2} \vert \nabla u_i \vert^2
&&= \int_{K_\rho} \vert x-y \vert^{-N+2-\alpha/2}\int_0^t \vert \nabla u_i \vert^2
\nonumber \\
&&
 \le C(C_s,s,N,t,\Vert u_i^0\Vert_{L^2(K_{2\rho})},\Vert Q  \Vert_{L^{s/(s-1)}(K_{2\rho} \times (0,t))}) \rho^{\alpha/2} 
\label{5}
\end{eqnarray}
for any $y \in K_\rho$.

We now can think of appealing to Lemma 4.4. page 61 in \cite{Lady}: Suppose that a function $u \ge 0$ satisfies for all $y \in K_\rho$
$$ \int_{K_\rho} \vert x-y\vert^{-N+m-\alpha/2}u^m \le C \rho^{\alpha/2}$$
with $\alpha >0$ and $1 < m \le 2$. Then, for any $\zeta \in W^{1,m}(K_\rho)$ with zero trace on the boundary $\partial K_\rho$, the following inequality holds true:
$$\int_{K_\rho}  u^m\zeta^2 \le C_1(C,N,m,\alpha) \rho^{2\alpha/m} \int_{K_\rho} u^{m-2} \vert \nabla \zeta\vert^2.$$
Unfortunately, the latter result is proved using  several H\"older's inequalities and the argument cannot be directly transposed to our time-dependent framework.

%\medskip

\subsubsection{The stationary case}

%\noindent{\bf Step 3: the stationary case.}\\
We restrict for some lines the study to the stationary case.
The reader can check straightforward that our previous computations remain almost unchanged for the elliptic setting.
We now can apply the result of Lemma 4.4. page 61 in \cite{Lady} mentioned above.
It allows to infer from \eqref{5} with $m=2$ that
\begin{equation}
  \int_{K_\rho} \vert v_i \vert^2 \vert \nabla u_i \vert^2 \le C(N,C_s,s,\Vert Q  \Vert_{L^{s/(s-1)}(K_{2\rho} )}) \rho^{(s-N)/2s} 
 \int_{K_\rho} \vert \nabla v_i \vert^2 .
\label{6}
\end{equation}

\medskip

\noindent{\bf Conclusion for the stationnary case.}\\
Estimate \eqref{6} under assumptions (i)-(ii) allow the control of $J$, thus of $I$ in \eqref{0} provided that $\rho$ is small enough.
Indeed, by combining all the inequalities above, we obtain that 
\begin{eqnarray}
&&
(\delta_1-2\epsilon_1- \frac{K_{1,+}}{2\epsilon_1}C(N,C_s,s,\Vert Q_1  \Vert_{L^{s/(s-1)}(K_{2\rho} )}) \rho^{(s-N)/2s}  - \frac{\ell (K_{2,1}^+)^2 }{4\epsilon_4})\int_{K_\rho}|\nabla v_1|^2\, dxds
 \nonumber \\
&&\quad  
+(\delta_2-2\epsilon_2-  \frac{K_{2,+}}{2\epsilon_2}C(N,C_s,s,\Vert Q_2  \Vert_{L^{s/(s-1)}(K_{2\rho} )} )\rho^{(s-N)/2s} -\frac{\ell (K_{1,2}^+)^2 }{4\epsilon_3})\int_{K_\rho}|\nabla v_2|^2\, dxds  
\nonumber \\
&& \quad
+(K_{1,1}^- - \epsilon_3)\int_{K_\rho}T_\ell(u_1) |\nabla v_1|^2 \, dxds
+(K_{2,2 }^- - \epsilon_4)\int_{K_\rho}T_\ell(u_2) |\nabla v_2|^2  \, dxds
\nonumber \\
&&
\leq 0.
\label{ineq2}
\end{eqnarray}
Hence, assuming $s>N$ and $\rho$ small enough, we can conclude that $v_i=0$ almost everywhere. The local in space uniqueness is proved.
The result reads as follows.

\begin{Prop}[Stationary case, local in space uniqueness] \label{prop2}
%Assume $\Omega \subset \mathbb{R}^2$. 
Assume that two weak solutions $(u_1,u_2)$ and $(\bar u_1,\bar u_2)$ of the elliptic version of Problem \eqref{problem} coincide on $\partial K_\rho $, for an open sphere $K_\rho \subset \Omega$ of radius $\rho$.
Assume\footnote{This first assumption is also made in \cite{CRR}.} 
$$ \frac{(K_{1,2}^+)^2}{K_{1,1}^-} < \frac{4\delta_2}{\ell}, \quad
 \frac{(K_{2,1}^+)^2}{K_{2,2}^-} < \frac{4\delta_1}{\ell}.$$
 Assume that the source terms satisfy \eqref{2b} with $s>N$. 
Assume\footnote{See Prop. \ref{prop1}. This second assumption is weaker than the one in \cite{CRR} and may be obtained without computing $g(s)$ if $N=2$.} further that $K$ is such that $(u_1,u_2)$ and $(\bar u_1,\bar u_2)$ actually belong to $W^{1,s}(\Omega)$.
Then 
$$(u_1,u_2)=(\bar u_1,\bar u_2) \  \mbox{a.e. in }  K_\rho .$$
\end{Prop}
\smallskip

Notice that the two results issued from \cite{Lady} in the latter proof still hold true if $\partial K_\rho \cap \partial \Omega \ne \emptyset$. The interested reader may check easily that all the other computations remain true replacing $K_\rho$ by $\Omega$. 
In this case, \eqref{4} reads
\begin{equation}
\int_{\Omega} \vert \nabla u_i \vert^2 \le C \bigl(\delta_i- \frac{\ell {K_{i,-i}^+}^2}{4K_{i,i}^-} \bigr)^{-1}\Vert Q  \Vert_{L^{s/(s-1)}(\Omega )}
\label{4b}
\end{equation}
and \eqref{6} simplifies into:
\begin{equation}
 \int_{\Omega} \vert v_i \vert^2 \vert \nabla u_i \vert^2 \le C(C_s,s,\Vert Q  \Vert_{L^{s/(s-1)}(\Omega )}) \vert \mbox{diam}(\Omega)\vert^{(s-2)/s} 
 \int_{\Omega} \vert \nabla v_i \vert^2 .
\label{6b}
\end{equation}

It follows that the latter result may be viewed as a global uniqueness result in the whole domain $\Omega$ provided its diameter is sufficiently small and provided its boundary is sufficiently regular\footnote{Indeed, in that case, we have to bring the proof of Lemma 4.3 in \cite{Lady} from the sphere $K_\rho$ to the whole $\Omega$.}.

\begin{Prop}[Stationary case, global uniqueness in a small and smooth domain] \label{prop3}
Assume $\Omega$ is a smooth domain of  $\mathbb{R}^N$. 
Assume
$$ \frac{(K_{1,2}^+)^2}{K_{1,1}^-} < \frac{3\delta_2}{\ell}, \quad
 \frac{(K_{2,1}^+)^2}{K_{2,2}^-} < \frac{3\delta_1}{\ell}.$$
 Assume that the source terms satisfy \eqref{2b} with  $s>N$, $\rho = \mbox{diam}(\Omega)$ and $C$ sufficiently small with regard to $\delta_i$.
Assume further that $K$ satisfy the assumptions in Proposition \ref{prop1} for $s>N$.
Then  the weak solution of the elliptic version of Problem \ref{problem} is unique in $W^{1,s}(\Omega)$.
\end{Prop}

\subsubsection{The time-dependent setting}
For turning back to the general setting, we can assume an additional hypothesis for ensuring that 
\eqref{6} remains true. 
Reading the proof of Lemma 4.4. page 61 in \cite{Lady} shows that assuming further that
\begin{equation}
\nabla u_i,\ i=1,2, \mbox{ belongs to } (L^\infty(0,T;L^2(\Omega)))^N
\label{7}
\end{equation}
ensures that all the results presented in the latter subsection extend to the time dependent case provided that the assumptions made on the source terms $Q_i$ also hold true for the initial data 
$\Vert u_i^0\Vert_{L^2}$ (see estimate \eqref{5}).

\bigskip

%%%%%%%%%%%%%%%%%%%%%
\section{Enhanced regularity and maximum principle}
%%%%%%%%%%%%%%%%%%%%%

In \cite{CRR}, the authors consider the question of the boundedness of the solutions of \eqref{problem}--\eqref{initial} without using their enhanced regularity result: assuming solely that the assumptions in Theorem \ref{Theo1} fulfilled, they prove that there exists source terms $Q_i \in L^2(0,T;(H^1(\Omega)'))$, $i=1,2$, such that the system \eqref{model} completed by the initial and boundary conditions \eqref{initial} admits a weak global solution such that, for  any $ T > 0 $, $(u_i-u_{i,D})_{i=1,2} \in W(0,T)^2$ and  the following maximum principle holds true:
$$ 0\leq u_i(t,x)  \leq \ell \quad \textrm{for a.e. } x \in \Omega , \textrm{ for all } t \in (0,T) \textrm{ and for all }   i=1,2.$$

In the present section, we aim at exploring if the enhanced regularity obtained in Proposition \ref{prop1} 
may be exploited for stating a maximum principle holding for a class of source terms.
It turns out that such a result  holds true provided that the regularity enhancement is sufficient (and actually quite important, see below).

We state and prove the following result.

\begin{Prop}[Explicit bound of the solutions of \eqref{problem}--\eqref{initial}]
\label{prop4}
Let $\ell^0 >0$. 
Assume that the source terms are such that $Q_i=Q_i(t,x,u_i)$ with $Q_i(t,x,y) \in L^s(0,T;W^{-1,s}(\Omega))$ a.e. $y \in \mathbb{R}$, $Q_i(t,x,y)\ge 0$ if $y \le 0$ and $Q_i(t,x,y)\le 0$ if $y \ge \ell^0$, a.e. in $\Omega_T$.
Assume that the initial and boundary data are such that
$$ 0 \le u_i^0 \le \ell^0 \mbox{ a.e. in }\Omega, \quad 
0 \le u_i^D \le \ell^0 \mbox{ a.e. in }(0,T).
$$
Assume that the assumptions in Proposition \ref{prop1} hold true with \footnote{Bear in mind that this specification for the regularity characteristic $s$  may be specified  using the assumptions of Proposition 1 in \cite{CRR}.} $s = 2N/(N-1)$. 
For any  $m>1$, there exists\footnote{See its explicit value in \eqref{est5} below.} a real number $C=C(\ell^0,m,s,K_{ij},\delta_i,\ell,\Vert u_i^0\Vert_{W^{1,s}(\Omega)},N)$ such that, if $T\vert \Omega\vert \le C$ then
$$0 \le u_i(t,x) \le m\ell^0 \mbox{ a.e. in } \Omega_T, \  i=1,2.$$
\end{Prop}
\begin{Rem}
The bound in Proposition \ref{prop4} depends in particular on $T$, on $\vert \Omega\vert$ and on $\Vert u_i^0\Vert_\infty$.  
Such a dependence is classical for quasilinear parabolic equations (the interested reader may for instance check that the result in Proposition \ref{prop4} is slightly better than those in Section 6, Chapter 2, in \cite{LadyP} ; notice that Zhou obtained in Theorem 1 in \cite{Zhou} a better result without any condition on $\vert\Omega\vert T$, but only for classical solutions of a nonlinear parabolic equation).
It means that if the quantities $T$,  $\vert \Omega\vert$ and  $\Vert u_i^0\Vert_\infty$ are sufficiently small (especially $\Vert u^0 \Vert_\infty \le \ell^0 < \ell$) we have
$$ 0 \leq u_i(t,x)  \leq \ell \quad \textrm{for a.e. } x \in \Omega , \textrm{ for all } t \in (0,T) \textrm{ and for all }   i=1,2$$
and that the system \eqref{problem} is actually the system \eqref{model}.
\end{Rem}
\begin{Rem}
For $N=2$, the result in Proposition \ref{prop4}  holds true provided that the solutions belong to $L^s(0,T;W^{1,s}(\Omega))$ with $s>4$. It means that the solutions are H\"older continuous in space.
Such a regularity is just above the one assumed to prove the uniqueness result in \cite{CRR}.
\end{Rem}

\bigskip

\begin{proof}
The inequality $0 \le u_i(t,x)$ almost everywhere in $\Omega_T$ was already obtained in Theorem \ref{Theo1}. 
Let $k \ge \ell^0$.
We write the variational formulation of the equation
$$\partial_t u_i - \mathrm{div} \bigl( (\delta_i + K_{ii}T_\ell(u_i))\nabla u_i \bigr) - \mathrm{div}(K_{i -i}T_\ell(u_i) \nabla u_{-i}) = Q_i$$
for the test function $(u_i -k)^+=\max\{ 0, u_i -k\}$.
Integrating by parts we get
\begin{eqnarray*}
&& \frac{1}{2} \Vert (u_i - k )^+ \Vert_{L^\infty(0,T;L^2(\Omega))}^2 + (\delta_i + K_{ii}^- \min\{k,\ell\}) \Vert \nabla (u_i -k)^+ \Vert_{(L^2(0,T;L^2(\Omega)))^N}^2
\\
&& \qquad \qquad
\le \Bigl\vert \int_{\Omega_T} K_{i -i} T_\ell(u_i) \nabla u_{-i} \cdot \nabla (u_i-k)^+ \Bigr\vert
\le K_{i -i}^+ \min\{k,\ell\} M_s^2 \mu_i(k)^{(s-2)/s}
\end{eqnarray*}
where we set
$$ \mu_i(k) = \mbox{mes}\{ (t,x)\in \Omega_T \mbox{ s.t. } u_i(t,x) >k\} = \int_0^T \int_\Omega \chi_{\{u_i>k\}}(t,x)\, dxdt$$
and $M_s$ is the real number such that
$$ \Vert  u_i \Vert_{L^s(0,T;W^{1,s}(\Omega))} \le M_s.$$
Notice that, according to the computations in \cite{CRR}, the dependence
$$M_S=M_s(s,K_{ij},\delta_i, T,\ell,\Vert u_i^0\Vert_{W^{1,s}(\Omega)})$$
is explicit.
It follows that
\begin{gather*}
 \Vert (u_i - k )^+ \Vert_{L^\infty(0,T;L^2(\Omega))}^2 + \Vert \nabla (u_i -k)^+ \Vert_{(L^2(0,T;L^2(\Omega)))^N}^2 \le C \mu_i(k)^{(s-2)/s}, \\
C = \frac{2}{\min\{ 1, \delta_i + K_{ii}^- \min\{k,\ell\}\}} K_{i -i}^+ \min\{k,\ell\} M_s^2 
\end{gather*}
and
\begin{gather}
 \Vert (u_i - k )^+ \Vert_{L^\infty(0,T;L^2(\Omega))} + \Vert \nabla (u_i -k)^+ \Vert_{(L^2(0,T;L^2(\Omega)))^N} \le C_i \mu_i(k)^{(s-2)/2s}, 
 \nonumber \\
C_i= \frac{\sqrt{2K_{i -i}^+ \min\{k,\ell\} }M_s}{\min\{ 1, \sqrt{\delta_i + K_{ii}^- \min\{k,\ell\}}\}}.
\label{estmu}
\end{gather}

We now aim at exploiting \eqref{estmu}, noticing especially that $C_i$ does not depend on $\ell^0$.
Let $m>1$, $m'>0$.
 Let $k_n=m\ell^0(1+m'-2^{-n})$ for  any $n \in \mathcal{N}$, $\mathcal{N}=\{n \in \mathbb{N} \mbox{ such that } k_n \ge m\ell^0\}$.
 Let $n_0=\min \mathcal{N}$.

First, using classical Sobolev injections, we notice that if $q \in [2,2N/(N-2)]$ and $r> 2$ are such that
$$ \frac{1}{r}+\frac{N}{2q} = \frac{1}{2}$$
then there exists some $\alpha>0$ such that the following interpolation inequality holds true:
$$ \Vert u_i \Vert_{L^r(0,T;L^q(\Omega))} \le \alpha \Vert u_i\Vert_{L^\infty(0,T;L^2(\Omega))}^{1-2/r} \Vert \nabla u_i \Vert_{(L^2(\Omega_T))^N}^{2/r}$$
and then, according to the Young inequality, there exists some $\beta>0$ such that
$$ \Vert u_i \Vert_{L^r(0,T;L^q(\Omega))} \le \beta \bigl( \Vert u_i\Vert_{L^\infty(0,T;L^2(\Omega))} + \Vert \nabla u_i \Vert_{(L^2(\Omega_T))^N}\bigr) .$$
Thus it follows from \eqref{estmu} that
\begin{equation}
\Vert (u_i - k_n)^+ \Vert_{L^r(0,T;L^q(\Omega))} \le C_i\beta \mu_i(k_n)^{(s-2)/2s}.
\label{est1}
\end{equation}
Next, we fix $q=r=N+2$. Since 
$$ (u-k_n)^+ \chi_{\{ (u-k_{n+1})^+\}} =  (u-k_n) \chi_{\{ (u-k_{n+1})^+\}} \ge  (k_{n+1}-k_n) \chi_{\{ (u-k_{n+1})^+\}} ,$$
we have
$$ \int_{\Omega_T} \vert (u_i-k_n)^+ \vert^r \chi_{\{ (u_i-k_{n+1})^+\ne 0\}} \, dxdt \ge (k_{n+1}-k_n)^r \mu_i(k_{n+1})$$
and thus
\begin{equation}
(k_{n+1}-k_n) \mu_i(k_{n+1})^{1/r} \le \Vert (u_i-k_n)^+ \Vert_{L^r(0,T;L^r(\Omega))}
\label{est2}
\end{equation}
where $k_{n+1}-k_n=m\ell^0/2^{n+1}$.
We infer from \eqref{est1}--\eqref{est2} that the sequence $(v_n)_{n \in \mathcal{N}}=(\mu_i(k_n))_{n \in \mathcal{N}}$ is such that
\begin{equation}
v_{n+1} \le (2^{n+1}C_i \beta v_n^{(s-2)/2s} /m{\ell^0})^r.
\label{est3}
\end{equation}
One may check that such a sequence satisfies $v_n \to 0$ as $n \to \infty$ if the two following conditions hold true:
$$\left\{ \begin{array}{l}
r(s-2)/2s :=1+\zeta >1 ,\\
{\displaystyle v_{n_0} \le  (m{\ell^0})^{r/\zeta} 2^{-r(1/\zeta+1/\zeta^2)}(C_i\beta)^{-r/\zeta} .}
\end{array}\right.$$
The first condition is reached as soon as $s>2N/(N-1)$.
The second condition is ensured if
\begin{equation}
 v_{n_0} = \mu_i(k_{n_0}) \le  (m{\ell^0})^{r/\zeta} 2^{-r(1/\zeta+1/\zeta^2)}(C_i\beta)^{-r/\zeta} .
 \label{est4}
 \end{equation}
Replacing $k_n$ by $ \ell^0$ and $k_{n+1}$ by $k_{n_0}$ in \eqref{est1}--\eqref{est2}, we get
 $$ (m-1)\ell^0 \mu_i(k_{n_0})^{1/r} \le C_i \beta \mu_i(\ell^0)^{(1+\zeta)/r} $$
and, since $\mu_i(\ell^0) \le T\vert \Omega\vert$,
$$ \mu_i(k_{n_0})^{1/r} \le \frac{C_i\beta}{(m-1)\ell^0} T^{(1+\zeta)/r}\vert \Omega\vert^{(1+\zeta)/r} .$$
Hence the condition \eqref{est4} is ensured if
\begin{gather}
 \frac{(C_i\beta)^r}{(m-1)^r(\ell^0)^r} T^{1+\zeta}\vert \Omega\vert^{1+\zeta} \le  (m{\ell^0})^{r/\zeta} 2^{-r(1/\zeta+1/\zeta^2)}(C_i\beta)^{r/\zeta} 
 \nonumber \\
 \Leftrightarrow T^{1+\zeta}\vert \Omega\vert^{1+\zeta} \le (m-1)^r m^{r/\zeta} 2^{-r(1/\zeta+1/\zeta^2)}(C_i \beta)^{r(1/\zeta-1)} (\ell^0)^{r(1+1/\zeta)} .
 \label{est5}
 \end{gather}
 If \eqref{est5} is satisfied, passing to the limit $n \to \infty$, we obtain
 $$ \lim_{n \to \infty} \mu_i(k_n) = \mu_i((1+m')m\ell^0) = 0$$
 for any $m'>0$, thus $ \mu_i(m\ell^0) =0$ that is $u_i(t,x) \le m\ell^0$ a.e. in $\Omega_T$.
\end{proof}

\section{Concept of confined solution}
\label{confined}

The result in Proposition \ref{prop4} has two weaknesses: a quite large regularity enhancement and a limitation of the size of the domain of interest $\Omega_T$ are necessary.
Whatever, it is proved in \cite{CRR} that these assumptions are not necessary, at least for a source term: there exists a confined solution of the problem \eqref{model}, \eqref{initial} in the following sense.

\textcolor{blue}{
\begin{Def}
The problem \eqref{model} completed by appropriate boundary and initial conditions admits a confined solution if there exists a source term $Q \in (L^2(0,T;(H^1(\Omega))'))^m$ and $u \in (W(0,T))^m$ such that $u_i$ solves
$$ \partial_t u_i - \nabla \cdot \big(\delta_i  \nabla u_i + u_i \sum_{j=1}^m K_{i,j}  \nabla u_j \big)  = Q_i \mbox{ in } \Omega_T$$
and $u_i$  is bounded almost everywhere in $\Omega_T$, $i=1,..,m$.
\end{Def}
}

\textcolor{blue}{
The advantage of this definition is that the term `confined' clearly corresponds to the construction of the solution which is forced  to remain  bounded by the penalization method. 
Another asset is that it sometimes corresponds to a physical interpretation of the confinement.
}

This latter point requires some precisions. 
In \cite{CRR}, the physical interpretation is detailed for the example of aquifer modelling.

\begin{figure}[tbp]
\centering \includegraphics[scale=0.42,angle=0]{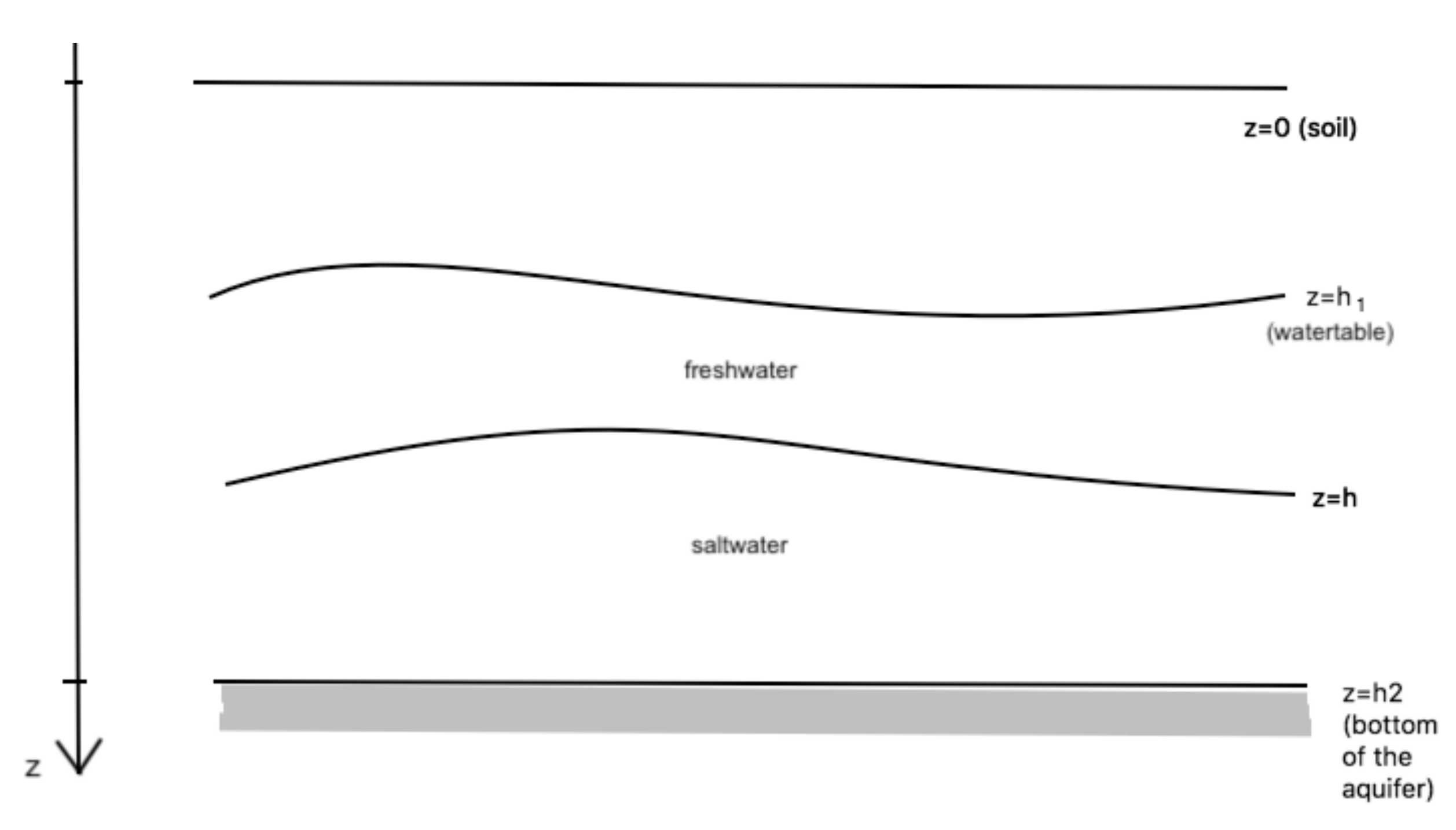}
% \centering \includegraphics[scale=0.7,angle=0]{figures/transition-zone}
  \caption{Aquifers modeling}
  \label{fig1}
\end{figure}

Define the depths $h$, $h_1$ and $h_2$ as in Figure \ref{fig1}.
The saltwater intrusion in the aquifer may be modeled by the following system (see \cite{CDR1}):
\begin{eqnarray}
\partial_t h - \delta \Delta h + \alpha \nabla \cdot \bigl( (h_2-h) \nabla h) \bigr) - \nabla \cdot \bigl( (1-\alpha)(h_2-h) \nabla h_1 \bigr) = 0,
\label{IntDep1}
\\
 \partial_t h_1 - \delta \Delta h_1 - \nabla \cdot \bigl( (1-\alpha)(h_2-h_1) \nabla h_1) \bigr) - \alpha \nabla \cdot \bigl( (h_2-h) \nabla h \bigr) = 0.
\label{IntDep2}
\end{eqnarray}
Complete the latter system by  initial and Dirichlet boundary conditions. 
Set  $u_1=h-h_1$ and $u_2=h_2-h$.
The system \eqref{IntDep1}--\eqref{IntDep2} enters the formalism of \eqref{model}, \eqref{initial}.
Hence,  assuming the necessary conditions  for Theorem \ref{Theo1}, namely $\ell =h_2$ and $1- 4\delta/h_2 < \alpha \le 1$, we can prove the existence of a weak solution $u=(u_1,u_2)$, with nonnegative components, and thus of $h$ and $h_1$ solving \eqref{IntDep1}--\eqref{IntDep2} 
in any given space-time domain $\Omega_T$.
Assuming moreover that the  initial and Dirichlet boundary conditions respect the physical hierarchy of interface depths, $h_1 \le h \le h_2$ a.e. in $\Omega_T$, we prove in \cite{CRR} that there exists  a confined solution of this problem.
To this aim, since the physical intuition consists in trying to prove that $0 \le h_1$, that is $u_1+u_2 \le h_2$ a.e. in $\Omega_T$, we add an {\it ad hoc}  penalization term in the equation characterizing $s=u_1+u_2$, namely
\begin{eqnarray}
&& \partial_t s^\epsilon - \delta \Delta s^\epsilon 
- \nabla \Bigl( \bigl( U_0(s^\epsilon- u_1^\epsilon) + (1-\alpha) U_0(u_1^\epsilon) \bigr)  \nabla s^\epsilon \Bigr) 
\nonumber \\
&& \qquad\qquad
- \alpha \nabla \cdot \bigr( U_0(u_1^\epsilon-s^\epsilon)\nabla u_1^\epsilon)
- \epsilon^{-1} \nabla \cdot \bigl( U_0(s^\epsilon-u_1^\epsilon)\nabla U_0(s^\epsilon-h_2) \bigr)
=0,
\label{IntDep22}
\end{eqnarray}
where  $ U_0 (x) = \max(0,x)$, and we let $\epsilon \to 0$.

\textcolor{blue}{
The interesting point is that there exists a physical interpretation of the latter penalization process.
With the penalization term in \eqref{IntDep22}, we assume that the aquifer is highly permeable above the depth $z=0$, thus the very high averaged permeability, namely equal to $\epsilon^{-1}$, when the thickness $u_1+u_2$ of the water exceeds $h_2$.
At the first order, this very conductive layer acts like a confining layer, as emphasized by the bound} $u_1+u_2 \ge 0$ \textcolor{blue}{at the limit $\epsilon \to 0$.
}
The situation is comparable to the presence of a highly conductive layer, a shallow substratum, at the top of the aquifer, which acts as a drain, and where the flow has a predominantly horizontal direction (see \cite{Van}, \cite{Ritz}).
\textcolor{blue}{
The mathematically confined solution $(h_1,h)$ of \eqref{IntDep1}-\eqref{IntDep2} with $0 \le h_1 \le h \le h_2$ a.e. in $\Omega_T$,  appears as the weak solution of 
\begin{eqnarray}
&& \partial_t h - \delta \Delta h + \alpha \nabla \cdot \bigl( (h_2-h) \nabla h) \bigr) - \nabla \cdot \bigl( (1-\alpha)(h_2-h) \nabla h_1 \bigr)
%\nonumber \\
%&& \qquad \qquad
-\nabla \cdot \mathcal{Q}  = 0,
\label{IntF1}
\\
&& \partial_t h_1 - \delta \Delta h_1 - \nabla \cdot \bigl( (1-\alpha)(h_2-h_1) \nabla h_1) \bigr) - \alpha \nabla \cdot \bigl( (h_2-h) \nabla h \bigr) 
%\nonumber \\
%&& \qquad \qquad
-   \nabla \cdot \mathcal{Q} 
= 0,
\label{IntF2}
\end{eqnarray}
in $\Omega_T$ completed by initial and Dirichlet boundary conditions, where $\mathcal{Q}\in (L^2(\Omega_T))^N$ is such that
$$ h_1 \mathcal{Q} = 0 \mbox{ a.e. in } \Omega_T.$$
}

\begin{figure}[tbp]
\centering \includegraphics[scale=0.264,angle=0]{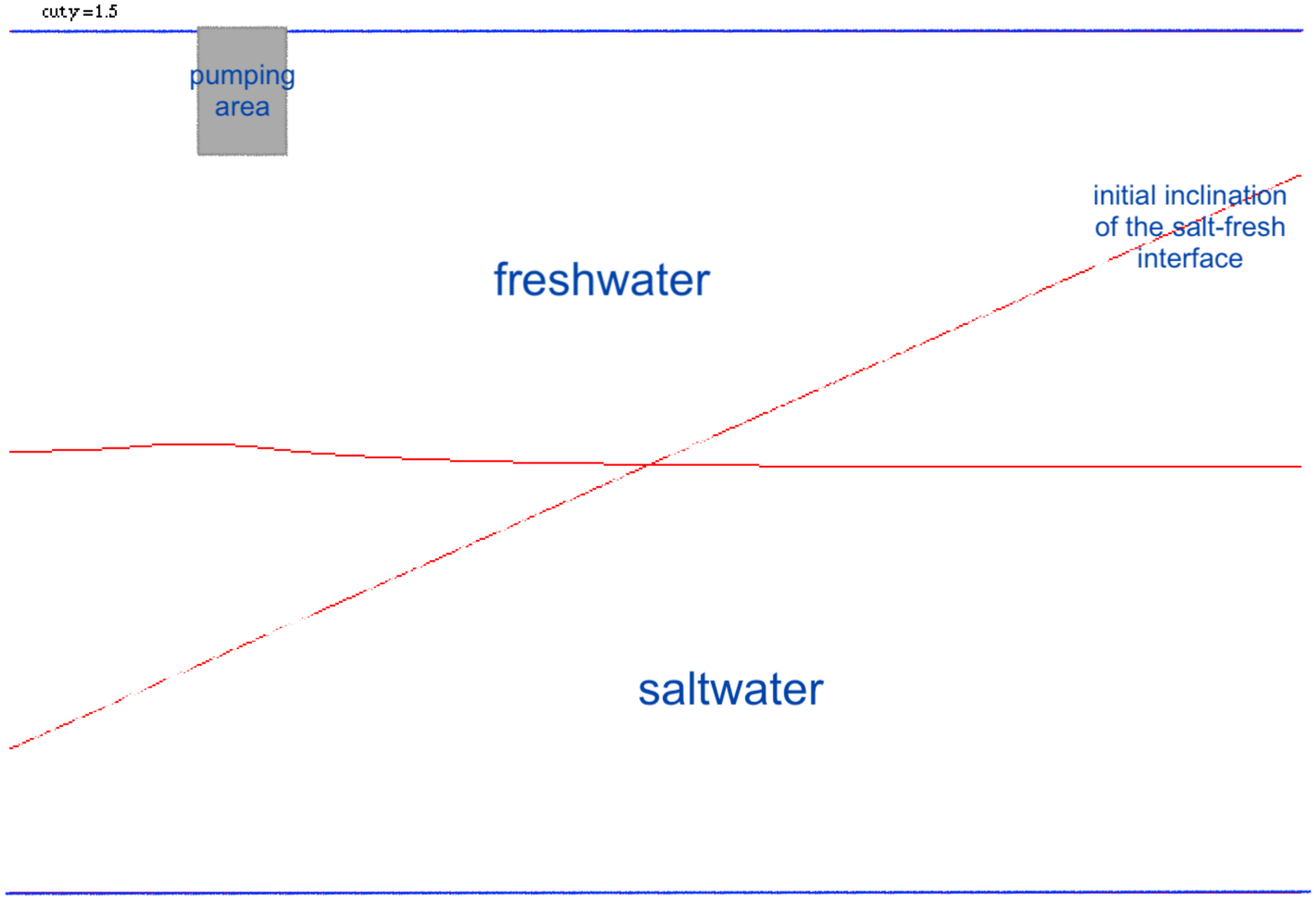}
\includegraphics[scale=0.26,angle=0]{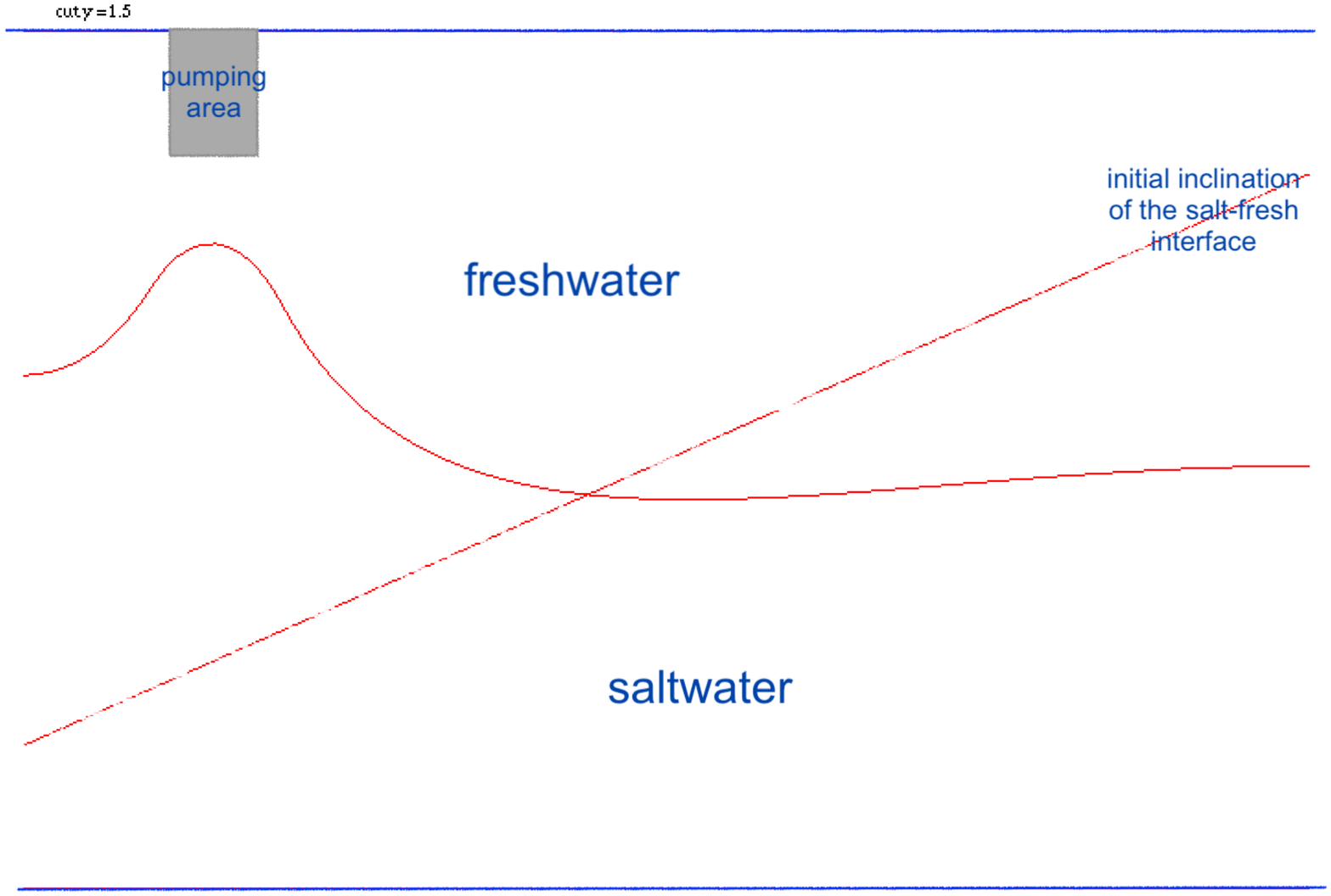}
  \caption{Keulegan experiment with pumping: confined aquifer (left) {\it versus} confined solution in an unconfined aquifer, {\it i.e.} solution confined  by an infinitely conductive upper layer (right). In the Keulegan experiment \cite{Keul}, the interface between salt-and freshwater is initially artificially inclined. Then  the interface should freely evolve due to the density contrast and the gravity effects  until horizontal stabilization. Here a pumping source term is added, thus the existence of a saltwater dome at the end of the computations.
The computations are done with the density contrast corresponding to seawater compared to clear water, $\alpha=  0.025$ and the same pumping rate.
}
  \label{f2}
\end{figure}

We would like to add an important note to avoid any confusion. 
Indeed, we have illustrated the concept of `confined solutions' by taking the example of aquifer models. Unfortunately, the term confinement is already used by hydrogeologists in the study of aquifers, but with a different meaning: in hydrogeology, a confined aquifer means that the reservoir is physically confined by an impermeable layer at its top and that it is fully saturated (that is $h_1=0$ here).
The mathematical model for the evolution of the salt interface $h$ and the hydraulic head $\Phi$ \emph{in a confined aquifer} is (see  \cite{CLR1})
\begin{eqnarray*}
&& \partial_t h - \delta \Delta h + \alpha \nabla \cdot \bigl( (h_2-h) \nabla h) \bigr) - \nabla \cdot \bigl( (1-\alpha)(h_2-h) \nabla \Phi \bigr)
= 0,
\\
&&  - \nabla \cdot \bigl( (1-\alpha)(h_2-h_1) \nabla \Phi) \bigr) - \alpha \nabla \cdot \bigl( (h_2-h) \nabla h \bigr) 
)
= 0.
\end{eqnarray*}
On the other hand, if we focus on the behaviour of \eqref{IntF1}-\eqref{IntF2} in a measurable subdomain where $h_1=0$ and $h_2-h \ge a_-$ for some $a_->0$, we notice that we can write $\mathcal{Q}=(1-\alpha)(h_2-h)  \nabla P$ with $P\in L^2(0,T;H_0^1(\Omega))$.
Hence, the \emph{confined solution of the unconfined aquifer model} solves: 
\begin{eqnarray*}
&& \partial_t h - \delta \Delta h + \alpha \nabla \cdot \bigl( (h_2-h) \nabla h) \bigr) - \nabla \cdot \bigl( (1-\alpha)(h_2-h) \nabla P \bigr)
= 0,
\\
&&  - \nabla \cdot \bigl( (1-\alpha)(h_2-h_1) \nabla P) \bigr) - \alpha \nabla \cdot \bigl( (h_2-h) \nabla h \bigr) 
%\\
%&& \qquad \qquad
- \nabla \cdot \bigl( (1-\alpha)(h_1-h) \nabla P) \bigr)
)
= 0
\end{eqnarray*}
with $P\in L^2(0,T;H_0^1(\Omega))$.
Simple numerical simulations show that the  two latter systems produce very different solutions (see {\it e.g.} Figure \ref{f2}).
However, in both models, the solutions remain confined (bounded), by an impermeable layer or by an infinitely permeable layer.

\bigskip

%%%%%%%%%%%%%%%%%%%%%%%%%%%%%%%%%%%%%%%%%%%%%%%%%
%%%%%%%%%%%%%%%%%%%%%%%%%%%%%%%%%%%%%%%%%%%%%%%
%\section*{References}

%%%%%%%%%%%%%%%%%%%%%%%%%%%%%%%%%%%%%%%%
\end{document}